\documentclass[11pt]{article}
\usepackage{amssymb,amsmath,amsthm}
\newtheorem{theorem}{Theorem}

\newtheorem{lemma}{Lemma}
\newtheorem{proposition}{Proposition}
\newtheorem{corollary}{Corollary}
\newtheorem{remark}{Remark}
\date{}
\numberwithin{equation}{section} \numberwithin{theorem}{section}
\numberwithin{lemma}{section} \numberwithin{corollary}{section}
\numberwithin{remark}{section} \numberwithin{proposition}{section}
\numberwithin{definition}{section}
\usepackage{color}
\usepackage{curves}
\usepackage[dvips]{graphics}
\usepackage{graphicx}
\usepackage{amsfonts}
\usepackage{amsmath}
\usepackage{amssymb}
\usepackage[normalem]{ulem}
\usepackage{t1enc}
\usepackage{fancybox}
\usepackage{epsfig}
\usepackage{minitoc}
\usepackage{fancyhdr}

\begin{document}

\newcommand{\n}{\noindent}

\newcommand{\vs}{\vskip}

\title{Uniqueness of solution for  a nonlinear  heterogeneous evolution dam problem}

\vs 0.5cm
\author{Messaouda Ben Attia$^{1}$,  Elmehdi Zaouche$^{2}$ and Mahmoud Bousselsal$^{3}$   \\
$^{1}$Universit\'{e} Kasdi Merbah\\
Laboratoire de Math\'{e}matiques Appliqu\'{e}es\\
 BP 511, Ouargla 30000, Algeria\\
benattiamessaouda1402@gmail.com\\
 $^{2}$University of El Oued\\
B. P. 789 El Oued 39000, Algeria\\
elmehdi-zaouche@univ-eloued.dz\\
elmehdizaouche45@gmail.com\\
$^{3}$Ecole Normale Sup\'{e}rieure\\
16050, Vieux-Kouba, Algiers, Algeria.\\
bousselsal55@gmail.com} \maketitle
\begin{abstract}
By choosing convenient test functions and using the method of
doubling variables, we prove the uniqueness of the solution to a
nonlinear  evolution dam problem in an arbitrary heterogeneous
porous medium of $\mathbb{R}^n (n\geq2)$  with an impermeable
horizontal bottom.
\end{abstract}

\par {\small {\bf Key words.} Test function; method of doubling
variables;  nonlinear
 evolution dam problem; heterogeneous porous medium; uniqueness. \\
\noindent {\bf  (2010) Subject Classifications}: 35A02, 35B35,
76S05.}

\section{Introduction}\label{s1}

Without loss of generality, we can assume that $n=2$. Let $\Omega$
be a bounded domain in $\mathbb{R}^{2}$ with horizontal bottom and
locally Lipschitz boundary $\partial\Omega:=\Gamma$ which represents
a porous medium and let $x=(x_1,x_2)$ be the generic point of
$\Omega$. Let $A$, $B$ and $D$ be real numbers such that $B>A$. The
boundary $\Gamma$ is divided into two parts such that one part
$\Gamma_1=[A,B]\times\{D\}$ is the impervious part which represents
the bottom of the dam and the other $\Gamma_2$ is the pervious part
which is a nonempty relatively open subset of $\Gamma$.  For a
positive real number $T,$ let $Q=\Omega\times(0,T)$ be the
space-time cylinder and let $\Sigma$ be the parabolic boundary
defined by $\Sigma_1\cup \Sigma_2$ where
$\Sigma_1=\Gamma_1\times(0,T)$ and
$\Sigma_2=\Gamma_2\times(0,T)=\Sigma_3\cup\Sigma_4$ with
$\Sigma_3=(\Gamma_2\times(0,T))\cap \{\varphi>0\}$,
$\Sigma_4=(\Gamma_2\times(0,T))\cap \{\varphi=0\}$ and $\varphi \in
C_{x}^{0,1}\cap C_{t}^{1}$ is a nonnegative function defined in
$\overline{Q}$ which represents the assigned pressure on
$\Gamma_2\times(0,T).$
 Let $a: \mathbb{R}\rightarrow \mathbb{R}$ be
a continuous function satisfying for some positive constants
$\lambda, \,\Lambda$ and $p>1,$
\begin{gather}\label{e1.1}
\forall r\in \mathbb{R}:  \quad \lambda |r|^{p}\leq a(r)r,\\
\forall r\in \mathbb{R}:   \quad |a(r)|\leq
\Lambda|r|^{p-1},\label{e1.2}\\
\forall r_1,r_2\in \mathbb{R}, \, r_1\neq r_2:\quad
(a(r_1)-a(r_2))(r_1-r_2)>0\label{e1.3}
\end{gather}
and let $h: (A,B)\rightarrow \mathbb{R}$ be a Lipschitz continuous
function of the variable $x_1$  such that  for two positive
constants $\underline{h}$ and $\overline{h},$
\begin{equation}\label{e1.1}
\underline{h}\leq h(x_1)\leq\overline{h} \quad\forall x_1\in (A,B).
\end{equation}
Moreover, let $g_0: \Omega\rightarrow\mathbb{R}$ be a measurable
function satisfying
\begin{equation}\label{e1.1}
0\leq g_0\leq1 \quad\text{ a.e. in }  \Omega.
\end{equation}
We set $\phi=\varphi+x_2$, and then we consider the following weak
formulation of a  nonlinear  heterogeneous evolution dam problem
associated with the initial data $g_0$:
\begin{equation*}{\bf (P)} \quad \left\{
\begin{aligned}
&\text{Find }
(u, g) \in L^p(0,T;W^{1,p}(\Omega))\times  L^\infty (Q) \text{ such that}:\\[0.2cm]
&\quad u \geq x_2, \; 0\leq g\leq 1,\; g(u-x_2) = 0  \quad \text{ a.e. in }Q, \\[0.1cm]
&\quad u=\phi  \hspace{0.34cm} \text{ on } \Sigma_2,\\
&\quad \displaystyle{\int_Q \big[h(x_1)(
a(u_{x_2})-ga(1))\xi_{x_2}+g\xi_t\big] \,dx dt+\int_\Omega g_0(x)\xi(x,0)\,dx}\leq0  \\
&\quad \forall \xi \in W^{1,p}(Q),\; \xi=0 \text{ on }
   \Sigma_3,\; \xi\geq0  \text{ on }
   \Sigma_4,\; \xi(x,T)=0 \,\text{ for a.e. }
x\in \Omega.
\end{aligned} \right.
\end{equation*}
In \cite{[Ly1]}, the author  established the existence of a solution
for the evolution dam problem related to an incompressible fluid
flow governed by a generalized nonlinear Darcy's law with Dirichlet
boundary conditions on some part of the boundary. He also proved in
\cite{[Ly2]} the continuity of solutions in $t$ for this problem.

\n For the homogeneous dam problem, the uniqueness of the solution
has been obtained in \cite{[C]} and \cite{[Ly1]} by the method of
doubling variables, respectively, for linear and generalized
nonlinear Darcy's laws. For a heterogeneous rectangular dam wet at
the bottom and dry near to the top, the uniqueness for a  linear
evolution dam problem has been proved in \cite{[LyZ]}, in both
incompressible and compressible flows, by an idea from \cite{[DF]}
in the homogeneous case. When $a(r)=r$, the uniqueness of the
problem $(P)$ in a rectangular porous medium has been obtained in
\cite{[Z1]} and \cite{[Z2]} by the method of doubling variables,
respectively, for incompressible and compressible flows.

\n In this paper, we choose convenient test functions and use the
method of doubling variables to prove the uniqueness of the solution
in a heterogeneous porous medium  of the evolution dam problem $(P)$
which is associated with an incompressible fluid governed by a
nonlinear Darcy's law. Our techniques are based on the uniqueness of
solutions obtained in \cite{[Ly1]} and \cite{[Z1]}. It should be
noted that our uniqueness result is new in the context of a
nonlinear evolution dam problem in an arbitrary heterogeneous
bounded domain of $\mathbb{R}^{n} (n\geq 2)$. In Section 2, we give
some properties of the solutions of $(P)$ and in Section 3, we state
and prove our uniqueness theorem that the solution of the problem
$(P)$ associated with the initial data $g_0$ is unique.

\section{Some properties of the solutions}\label{s2}

In this section, we will give some  properties of the solutions
which are useful in proving our main result.
\begin{lemma}(\cite{[Ly1]})\label{lma2.2}
Let $v\in W^{1,p}(Q)$ and $F\in W^{1,\infty}_{loc}(\mathbb{R}^2)$ be
functions satisfying
\begin{eqnarray}\label{e3.1}
&F(u-x_2, v) \in L^p(0,T; W^{1,p}(\Omega)), \; F(\phi-x_2,v)\in W^{1,p}(Q), \, F(z_1,z_2)\geq0 \, \text{ a.e. } (z_1,z_2)\in \mathbb{R}^2, \nonumber\\
&\text{ and either } \frac{\partial F}{\partial z_1}(z_1,z_2)\geq0
\, \text{ a.e. } (z_1,z_2)\in \mathbb{R}^2 \, \text{ or } \,
\frac{\partial F}{\partial z_1}(z_1,z_2)\leq0 \, \text{ a.e. }
(z_1,z_2)\in \mathbb{R}^2.\nonumber
\end{eqnarray}
Then, if $(u,g)$ is  a solution of $(P)$ and $\xi \in
\mathcal{D}(\overline{\Omega}\times(0,T)),$ we have
\begin{eqnarray}\label{e3.1}
&&\int_{Q}h(x_1)(a(u_{x_2})-ga(1))(F(u-x_2,v)\xi)_{x_2}+g(F(0,v)\xi)_t\, dx dt \nonumber\\
&&=\int_Qh(x_1)(a(u_{x_2})-ga(1))(F(\phi-x_2,v)\xi)_{x_2}+g(F(\phi-x_2,v)\xi)_t\,
dx dt.\nonumber
\end{eqnarray}
 In particular, if  $F(\phi-x_2,v)\xi=0$ on $\Sigma_2,$
\begin{eqnarray}\label{e3.1}
&&\int_{Q}h(x_1)(a(u_{x_2})-ga(1))(F(u-x_2,v)\xi)_{x_2}+g(F(0,v)\xi)_t\,
dx dt. \nonumber
\end{eqnarray}
\end{lemma}
The following corollary is an immediate consequence of Lemma 2.1.
\begin{corollary}\label{prop2.1}
Let $\epsilon>0$ and $k\geq0$ be real numbers and  let $\xi\in
\mathcal{D}(\mathbb{R}^{2}\times(0,T))$ such that $\xi\geq0$ and
$\xi=0$ on $\Sigma_3.$ If $(u,g)$ is a solution of $(P)$, we have
$$\int_{Q}h(x_1)a(u_{x_2})(\min(\frac{(u-x_2-k)^+}{\epsilon},1)\xi)_{x_2}\, dxdt=0.$$
\end{corollary}
Let us set
\begin{eqnarray}\label{e3.1}
&&\sigma_1=\overline{\Sigma}_2\cap\Sigma_1=(\overline{\Gamma}_2\cap\Gamma_1)\times(0,T),\nonumber\\
&&\sigma_2=\overline{\Sigma}_3\cap\Sigma_4=(\Gamma_2\times(0,T))\cap\overline{\{\varphi>0\}}\cap
\{\varphi=0\}\nonumber
\end{eqnarray}
and let us assume throughout the rest of the paper that $\sigma_1$
and $\sigma_2$ are $(1,q)$ polar sets of $\overline{Q}$ (see
\cite{[Ad]}), where $q$ is the conjugate exponent of $p$. Since the
empty set is the only $(1,q)$ polar set of $\overline{Q}$ in the
case $p>3$, then,  we can consider that $p\leq 3.$

We use a regularization by convolution with respect to the variables
$x_2$ and $t$ to prove the following proposition.
\begin{proposition}\label{prop2.1}
Let $\lambda\in [0,1]$ and let $\xi\in
\mathcal{D}(\mathbb{R}^{2}\times(0,T))$ such that $\xi\geq0$ and
$\xi=0$ on $\Sigma_1\cup \Sigma_3.$ If $(u,g)$ is a solution of
$(P)$, we have
\begin{eqnarray}\label{e3.1}
&&\int_{Q}\big\{h(x_1)(a(u_{x_2})-a(1))\xi_{x_2}+(\lambda-g)^{+}(h(x_1)a(1)\xi_{x_2}-\xi_t)\big\}\,dxdt\leq0.
\end{eqnarray}
\end{proposition}
\n \emph{Proof.} We apply Corollary 2.1 for  $k=0$ to get
\begin{eqnarray}\label{e3.1}
&&\int_{Q}h(x_1)
a(u_{x_2})(\min(\frac{u-x_2}{\epsilon},1)\xi)_{x_2}\,dxdt=0.
\end{eqnarray}
On the other hand, we have
\begin{eqnarray}\label{e3.1}
&&\int_{Q}h(x_1)a(1)(\min(\frac{u-x_2}{\epsilon},1)\xi)_{x_2}\,dxdt=
0\end{eqnarray}
 since $\min(\frac{u-x_2}{\epsilon},1)\xi=0$ on
 $\Sigma$ and $(h(x_1))_{x_2}=0$ a.e. in
 $Q.$ Subtracting (2.3) from and (2.2), we get
\begin{eqnarray}\label{e3.1}
&&\int_{Q}h(x_1)(a(u_{x_2})-a(1))(\min(\frac{u-x_2}{\epsilon},1)\xi)_{x_2}\,dxdt=0\nonumber
\end{eqnarray}
which can be written as
\begin{eqnarray}\label{e3.1}
&&\frac{1}{\epsilon}\int_{Q\cap\{u-x_2<\epsilon\}}\xi
h(x_1)(a(u_{x_2})-a(1))(u_{x_2}-1)\,dx
dt\nonumber\\
&&+\int_{Q}\min(\frac{u-x_2}{\epsilon},1)h(x_1)(a(u_{x_2})-a(1))\xi_{x_2}\,dxdt=0.
\end{eqnarray}
By (1.3) and the fact that $\xi h(x_1)\geq0$ a.e. in $Q$, the first
integral of (2.4) is nonnegative, then,
\begin{eqnarray}\label{e3.1}
&&\int_{Q}\min(\frac{u-x_2}{\epsilon},1)h(x_1)(a(u_{x_2})-a(1))\xi_{x_2}\,dxdt\leq0.
\end{eqnarray}
Letting $\epsilon\rightarrow0$ in (2.5), we obtain
\begin{eqnarray}\label{e3.1}
&&\int_{Q}h(x_1)(a(u_{x_2})-a(1))\xi_{x_2}\,dxdt\leq0
\end{eqnarray}
and then (2.1) holds for $\lambda=0.$  Also, the inequality (2.1)
holds for $\lambda=1$ since $0\leq g\leq1$ a.e. in $Q$ and $\xi$ is
a test function for $(P),$
\begin{eqnarray}\label{e3.1}
&&\int_{Q}\big\{h(x_1)(a(u_{x_2})-ga(1))\xi_{x_2}+g\xi_t\big\}\,dxdt\leq0.
\end{eqnarray}

\n Now, we will prove (2.1) for $\lambda\in(0,1).$   Without loss of
generality, we can assume that
$d(\text{supp}(\xi),\Sigma_1\cup\Sigma_3):=\varepsilon_0>0.$ Let us
set $A=((\mathbb{R}^2\times(0,T))\backslash\Sigma_1)\cup
\Sigma_3\cup \sigma_2$ and
$$A_{\varepsilon_0}=\big\{(x,t)\in \mathbb{R}^2\times(0,T)/ d((x,t),\Sigma_1\cup\Sigma_3\cup\sigma_2)>\frac{\varepsilon_0}{2}\big\}.$$
We extend $u$ (resp. $g$) on  $A\backslash Q$ by $x_2$  (resp. 1)
and still denote by $u$ (resp. $g$) this function. Also, the
function $f$ can be extended to a Lipschitz function on
$\mathbb{R}$,  still denote by $h$. We use a regularization by
convolution for $a(u_{x_2})$ and $g$ with respect to the variables
$x_2$ and $t, \, (a(u_{x_2}))_\varepsilon=\rho_\varepsilon\ast
a(u_{x_2}), \,g_\varepsilon=\rho_\varepsilon\ast g$ where
$\varepsilon \in (0,\frac{\varepsilon_0}{2}), \, \rho_\varepsilon\in
\mathcal{D}(\mathbb{R}\times(0,T)),\,
\text{supp}(\rho_\varepsilon)\subset B(0,\varepsilon)$ is a
regularizing sequence. We can use Fubini's theorem to write
\begin{eqnarray}\label{e3.1}
&&\int_{A_{\varepsilon_0}}\big\{h(x_1)((a(u_{x_2}))_\varepsilon-g_\varepsilon a(1))\xi_{x_2}+g_\varepsilon\xi_t\big\}\,dxdt\nonumber\\
&&=\int_{A_{\varepsilon_0}}\Big\{\int_{\mathbb{R}\times(0,T)}(a(u_{x_2})-ga(1))(x_1,x_2-y,t-s)\rho_{\varepsilon}(y,s)\,dyds\Big\}h(x_1)\xi_{x_2}\,dx
dt\nonumber\\
&&+\int_{A_{\varepsilon_0}}\Big\{\int_{\mathbb{R}\times(0,T)}g(x_1,x_2-y,t-s)\rho_{\varepsilon}(y,s)\,dyds\Big\}\xi_{t}\,dx
dt \nonumber\\
&&=\int_{\mathbb{R}\times(0,T)}\rho_{\varepsilon}(y,s)\Big\{\int_{A_{\varepsilon_0}}h(x_1)(a(u_{x_2})-ga(1))(x^{\prime},x_2-y,t-s)\xi_{x_2}\,dxdt\Big\}\,dyds\nonumber\\
&&+\int_{\mathbb{R}\times(0,T)}\rho_{\varepsilon}(y,s)\Big\{\int_{A_{\varepsilon_0}}g(x_1,x_2-y,t-s)\xi_{t}\,dxdt\Big\}\,dyds.\nonumber
\end{eqnarray}
Then, if we make the change of variables  $z=x_2-y$ and $\tau=t-s,$
we get
\begin{eqnarray}\label{e3.1}
&&\int_{A_{\varepsilon_0}}\big\{h(x_1)((a(u_{x_2}))_\varepsilon-g_\varepsilon a(1))\xi_{x_2}+g_\varepsilon\xi_t\big\}\,dxdt\nonumber\\
&&=\int_{B(0,\varepsilon)}\rho_{\varepsilon}(y,s)\Big\{\int_{A_{\varepsilon_0}}h(x_1)(a(u_{z})-
ga(1))(x_1,z,t)(\xi(x_1,z+y,\tau+s))_{z}\,dx_1dzd\tau\Big\}\,dyds\nonumber\\
&&+\int_{B(0,\varepsilon)}\rho_{\varepsilon}(y,s)
\Big\{\int_{A_{\varepsilon_0}}g(x_1,z,t)(\xi(x_1,z+y,\tau+s))_{\tau}\,dx_1dzd\tau\Big\}\,dy ds \nonumber\\
&&=\int_{B(0,\varepsilon)}\rho_{\varepsilon}(y,s)\Big\{\int_{Q}h(x_1)(a(u_{z})-
ga(1))(x_1,z,t)(\xi(x_1,z+y,\tau+s))_{z}\,dx_1dzd\tau\Big\}\,dyds\nonumber\\
&&+\int_{B(0,\varepsilon)}\rho_{\varepsilon}(y,s)\Big\{\int_{Q}g(x_1,z,t)(\xi(x_1,z+y,\tau+s))_{\tau}\,dx_1dzd\tau\Big\}\,dy
ds. \nonumber
\end{eqnarray}
Observe that  $(x_1,z,\tau)\mapsto \xi(x_1,z+y,\tau+s)$ is a
nonnegative  function in $\mathcal{D}(\mathbb{R}^{2}\times(0,T))$
and vanishes on $\Sigma_1\cup \Sigma_3$ for all $(y,s)\in
B(0,\varepsilon).$ Therefore, since $\rho_\varepsilon\geq0,$ we
deduce from (2.7) that
\begin{eqnarray}\label{e3.1}
&&\int_{A_{\varepsilon_0}}\big\{h(x_1)((a(u_{x_2}))_\varepsilon-g_\varepsilon
a(1))\xi_{x_2}+g_\varepsilon\xi_t\big\}\,dxdt\leq0\nonumber
\end{eqnarray}
which can be written as
\begin{eqnarray}\label{e3.1}
&&\int_{A_{\varepsilon_0}}\big\{h(x_1)((a(u_{x_2}))_\varepsilon-a(1))\xi_{x_2}+(\lambda
-g_\varepsilon)(h(x_1)a(1)\xi_{x_2}-\xi_t)\big\}\,dxdt\leq0\nonumber
\end{eqnarray}
since
$$\int_{A_{\varepsilon_0}}h(x_1)a(1) \xi_{x_2}\, dxdt=\int_{A_{\varepsilon_0}}\lambda\xi_{t}\, dxdt=0.$$
Similarly, using (2.6), we arrive at
\begin{eqnarray}\label{e3.1}
&&\int_{A_{\varepsilon_0}}h(x_1)((a(u_{x_2}))_\varepsilon-a(1))\xi_{x_2}\,dxdt\leq0.\nonumber
\end{eqnarray}
Now, for any positive real number $\delta$, we set
$K_\delta=\min(\frac{(\lambda-g_\varepsilon)^+}{\delta},1)$  which
satisfies $K_\delta\in L^{p}_{\text{loc}}(A_{\varepsilon_0}), \,
K_{\delta x_2}, K_{\delta t}\in
L^{p}_{\text{loc}}(A_{\varepsilon_0}).$ By using the integration by
parts formula, we obtain
\begin{eqnarray}\label{e3.1}
&&\int_{A_{\varepsilon_0}}\big\{h(x_1)((a(u_{x_2}))_\varepsilon-a(1))\xi_{x_2}+K_\delta(\lambda
-g_\varepsilon)(h(x_1)a(1)\xi_{x_2}-\xi_t)\big\}\,dx dt\nonumber\\
&&-\frac{\delta}{2}\int_{A_{\varepsilon_0}}K_\delta^2(h(x_1)a(1)\xi_{x_2}-\xi_t)\,dx
dt \nonumber\\
&&=\int_{A_{\varepsilon_0}}\big\{h(x_1)((a(u_{x_2}))_\varepsilon-a(1))(K_\delta\xi)_{x_2}\nonumber\\
&&\quad\quad\quad\quad\quad+(\lambda
-g_\varepsilon)(h(x_1)a(1)(K_\delta\xi)_{x_2}-(K_\delta\xi)_t)\big\}\,dx dt\nonumber\\
&&+\int_{A_{\varepsilon_0}}h(x_1)((a(u_{x_2}))_\varepsilon-a(1))((1-K_\delta)\xi)_{x_2}\,
dx dt\nonumber
\end{eqnarray}
and since (2.6) and (2.7) remain valid,  respectively, for
$K_\delta\xi$ and $(1-K_\delta)\xi,$ it follows that
\begin{eqnarray}\label{e3.1}
&&\int_{A_{\varepsilon_0}}\big\{h(x_1)((a(u_{x_2}))_\varepsilon-a(1))\xi_{x_2}+K_\delta(\lambda
-g_\varepsilon)(h(x_1)a(1)\xi_{x_2}-\xi_t)\big\}\,dx dt\nonumber\\
&&-\frac{\delta}{2}\int_{A_{\varepsilon_0}}K_\delta^2(h(x_1)a(1)\xi_{x_2}-\xi_t)\,dx
dt\leq0.
\end{eqnarray}
Finally, we pass successively to the limit in (2.8) as
$\delta\rightarrow 0$  and then as $\varepsilon\rightarrow 0$ and
using Lebesgue's dominated convergence theorem, we obtain
\begin{eqnarray}\label{e3.1}
&&\int_{A_{\varepsilon_0}}\big\{h(x_1)(a(u_{x_2})-a(1))\xi_{x_2}+(\lambda
-g)^+(h(x_1)a(1)\xi_{x_2}-\xi_t)\big\}\,dx dt\leq0,\nonumber
\end{eqnarray}
and then (2.1) holds since $u=x_2,\, g=1$ a.e. in $A\backslash Q$
and $\varepsilon_0$ is  arbitrary.\qed

We use Corollary 2.1 and Proposition 2.1 and argue as in the proof
of [\cite{[Ly1]}, Lemma 5.2] to prove the following lemma:
\begin{lemma}\label{lma2.2}
Let $\chi$ be a function of $L^{\infty}(Q)$ satisfying
$$0\leq \chi\leq 1  \quad \text{ and } \quad h(x_1)a(1)\chi_{x_2}-\chi_t=0 \quad \text{ in } \mathcal{D}^{\prime}(Q).$$
Let $\xi, \xi_1, \xi_2\in \mathcal{D}(\mathbb{R}^2\times(0,T))$ such
that $\xi, \xi_1\geq0, \, \xi=\xi_1=0$ on $\Sigma_1\cup \Sigma_3, \,
\xi_2=0$ on $\partial Q$ and let $k, \lambda, \epsilon$ be
nonnegative real numbers such that $\epsilon>0$ and $\lambda\in
1-H(k)$ with $H$ denotes the maximal monotone graph associated to
the Heaviside function. Then, if $(u,g)$ is a solution of $(P),$ we
have
\begin{eqnarray}\label{e3.1}
&&\int_{Q}\Big\{h(x_1)a(u_{x_2})\big(\min\big(\frac{(u-x_2-k)^+}{\epsilon},1\big)\xi\big)_{x_2}+
(\lambda-g)^+(h(x_1)a(1)\xi_{1x_2}-\xi_{1t})\nonumber\\
&&+(\lambda-\chi)^+(h(x_1)a(1)\xi_{2x_2}-\xi_{2t})\Big\}\, dx dt\leq
C(u,k,\xi_1),\nonumber
\end{eqnarray}
where
\begin{eqnarray}\label{e3.1}
C(u,0,\xi_1)&=&-\int_{Q}h(x_1)(a(u_{x_2})-a(1))\xi_{1x_2}\, dx
dt \nonumber\\
&=&\lim_{\epsilon\rightarrow0}\int_{Q}h(x_1)(a(u_{x_2})-a(1))\big(\min\big(\frac{u-x_2}{\epsilon},1\big)\big)_{x_2}\xi_1\,
dx dt,\\
 C(u,k,\xi_1)&=& 0 \quad \text{ for } k>0. \nonumber
\end{eqnarray}
\end{lemma}

We use Lemma 2.1 and  employ the regularization by convolution with
respect to the variables $x_2$ and $t$ as in the proof of
Proposition 2.1,  we obtain by an argument similar to that in the
proof of [\cite{[Ly1]}, Lemma 5.3],
\begin{lemma}\label{lma2.2}
Let $\nu$ denote the outward unit normal to $\Gamma_1$ and let us
assume that  $(0,h(x_1)a(1)).\nu \leq 0$ on $\Gamma_1$. Let $\Psi$
be a function of $C^\infty(\mathbb{R})\cap C^{0,1}(\mathbb{R})$ such
that $\Psi(0)=0, \, \Psi^{\prime}\geq 0, \, \Psi\leq 1$ and let $k,
\lambda, \epsilon$ be the nonnegative real numbers defined in Lemma
2.2. Then, if $(u,g)$ is a solution of $(P)$ and $\xi\in
\mathcal{D}(\mathbb{R}^2\times(0,T)), \, \xi\geq 0, \,
(1-\Psi(u-x_2))\xi=0$ on $\Sigma_2,$ we have
\begin{eqnarray}\label{e3.1}
&&\int_{Q}\Big\{h(x_1)(a(u_{x_2})-\lambda
a(1))\big(\min\big(\frac{(k-(u-x_2))^+}{\epsilon},1\big)(1-\Psi(u-x_2))\xi\big)_{x_2}\nonumber\\
&&-(g-\lambda)^+(h(x_1)a(1)\xi_{x_2}-\xi_{t})\Big\}\, dx dt\geq
0.\nonumber
\end{eqnarray}

\end{lemma}

\section{Uniqueness of solution}\label{s2}
In this section we state and prove our uniqueness theorem.
\begin{theorem}
 Assume that (1.1)-(1.5) and that $(0,h(x_1)a(1)).\nu \leq 0$ on $\Gamma_1$  hold. Then,
the solution of the problem $(P)$ associated with the initial data
$g_0$ is unique.
\end{theorem}
We seek to obtain a comparison result  for solutions which allows us
to prove the uniqueness of the solution of the problem $(P)$. First,
we begin with the following two comparison lemmas of solutions.

\begin{lemma}
Let $B$ be a bounded open subset of $\mathbb{R}^{2}$ such that
either $B\cap\Gamma=\emptyset$  or $B\cap\Gamma$  is a Lipschitz
graph. For two solutions $(u_1,g_1)$ and $(u_2,g_2)$ to $(P),$  we
set $u_m=\min(u_1,u_2)$ and $g_M=\max(g_1,g_2).$ Then, for all $\xi
\in \mathcal{D}(B\times(0,T)), \, \xi\geq0, \,
\text{supp}(\xi)\cap(\Sigma_1\cup \Sigma_3)=\emptyset$ and for
$i=1,2,$ we have
\begin{eqnarray}\label{e3.1}
&&\int_{Q}\big\{h(x_1)(a(u_{ix_2})-a(u_{mx_2})-(g_i-g_M)a(1))\xi_{x_2}\nonumber\\
&&\quad\quad+(g_i-g_M)\xi_t\big\}\,dxdt\leq0.
\end{eqnarray}
\end{lemma}
\n \emph{Proof.} Let $(u_1,g_1)$ and $(u_2,g_2)$ be two solutions of
 $(P)$ and let $\xi$ be the function defined in Lemma 3.1. We
 define
 \begin{eqnarray}\label{e3.1}
&&\forall (x,t,y,s)\in \overline{Q\times Q}:\nonumber\\
&&\zeta(x,t,y,s)=\xi(\frac{x+y}{2},\frac{t+s}{2})
\rho_{1,\delta_1}(\frac{t-s}{2})\rho_{2,\delta_1}(\frac{x_1-y_1}{2})\rho_{3,\delta_2}(\frac{x_2-y_2}{2}),\nonumber
\end{eqnarray}
where $\delta_1, \delta_2$  are  positive real numbers,
$\rho_{1,\delta_1}, \,\rho_{2,\delta_1},\, \rho_{3,\delta_2}\in
\mathcal{D}(\mathbb{R}), \, \rho_{1,\delta_1}, \,
\rho_{2,\delta_1},\, \rho_{3,\delta_2}\geq0$ in $\mathbb{R}, \,
\int_{\mathbb{R}}\rho_{1,\delta_1}(t)\,dt=\int_{\mathbb{R}}\rho_{2,\delta_1}(t)\,dt=\int_{\mathbb{R}}\rho_{3,\delta_2}(t)\,dt=1,
\,\text{supp}(\rho_{1,\delta_1}),
\,\text{supp}(\rho_{2,\delta_1})\subset(-\delta_1,\delta_1),\,
\text{supp}(\rho_{3,\delta_2})\subset(-\delta_2,\delta_2)$ and
\begin{eqnarray}\label{e3.1}
&&\forall (x,y)\in (B\cap\Omega)\times (B\backslash\Omega), \quad
\rho_{2,\delta_1}(\frac{x_1-y_1}{2})\rho_{3,\delta_2}(\frac{x_2-y_2}{2})=0.\nonumber
\end{eqnarray}
Notice that, by choosing $\delta_1$ and $\delta_2$ small enough,
$\zeta\in \mathcal{D}(B\times(0,T)\times B\times(0,T))$ and
\begin{eqnarray}\label{e3.1}
&&\zeta=0 \quad \text{ on } \, ((\Sigma_1\cup \Sigma_3)\times Q)\cup
(Q\times \Sigma).
\end{eqnarray}
So,  applying Lemma 2.2 to $(u_1,g_1)$ with $k=u_2(y,s)-y_n, \,
\lambda=g_2(y,s), \,\xi(x,t)=\xi_1(x,t)=\zeta(x,t,y,s)$ and
$\xi_2(x,t)=0,$ we obtain for a.e. $(y,s)\in Q,$
\begin{eqnarray}\label{e3.1}
&&\int_{Q}\Big\{h(x_1)a(u_{1x_2})\big(\min\big(\frac{(u_1-x_2+y_2-u_2)^+}{\epsilon},1\big)\xi\big)_{x_2}\nonumber\\
&&+(g_2-g_1)^+(h(x_1)a(1)\zeta_{x_2}-\zeta_{t})\Big\}\, dx dt\leq
C(u_1,u_2-y_2,\zeta)\nonumber
\end{eqnarray}
and integrating over $Q,$ we get
\begin{eqnarray}\label{e3.1}
&&\int_{Q\times Q}\Big\{h(x_1)a(u_{1x_2})\big(\min\big(\frac{(u_1-x_2+y_2-u_2)^+}{\epsilon},1\big)\xi\big)_{x_2}\nonumber\\
&&+(g_2-g_1)^+(h(x_1)a(1)\zeta_{x_2}-\zeta_{t})\Big\}\, dx
dtdyds\leq \int_{Q}C(u_1,u_2-y_2,\zeta)\, dy ds.
\end{eqnarray}
On the other hand, applying Lemma 2.3 to $(u_2,g_2)$ with
$k=u_1(x,t)-x_2, \, \lambda=g_1(x,t), \,\xi(y,s)=\zeta(x,t,y,s)$ and
$\Psi=0,$ we have for a.e. $(x,t)\in Q,$
\begin{eqnarray}\label{e3.1}
&&\int_{Q}\Big\{h(y_1)(a(u_{2y_2})-g_1
a(1))\big(\min\big(\frac{(u_1-x_2-u_2+y_2)^+}{\epsilon},1\big)\zeta\big)_{y_2}\nonumber\\
&&-(g_2-g_1)^+(h(y_1)a(1)\zeta_{y_2}-\zeta_{t})\Big\}\, dy ds\geq 0.
\end{eqnarray}
Using (3.2) and the fact that the function $(y,s)\mapsto
h(y_1)g_1a(1)$ does not depend on $y_2$, we find
$$\int_{Q}h(y_1) g_1
a(1)\big(\min\big(\frac{(u_1-x_2-u_2+y_2)^+}{\epsilon},1\big)\zeta\big)_{y_2}\,
dyds=0,$$
 therefore, (3.4) can be written as
\begin{eqnarray}\label{e3.1}
&&\int_{Q}\Big\{h(y_1)a(u_{2y_2})\big(\min\big(\frac{(u_1-x_2-u_2+y_2)^+}{\epsilon},1\big)\zeta\big)_{y_2}\nonumber\\
&&-(g_2-g_1)^+(h(y_1)a(1)\zeta_{y_2}-\zeta_{t})\Big\}\, dy ds\geq
0.\nonumber
\end{eqnarray}
By integrating over $Q$, we obtain
\begin{eqnarray}\label{e3.1}
&&-\int_{Q\times Q}\Big\{h(y_1)a(u_{2y_2})\big(\min\big(\frac{(u_1-x_2-u_2+y_2)^+}{\epsilon},1\big)\zeta\big)_{y_2}\nonumber\\
&&-(g_2-g_1)^+(h(y_1)a(1)\zeta_{y_2}-\zeta_{t})\Big\}\, dx dt dy
ds\leq 0.
\end{eqnarray}
Since
$$\min\big(\frac{(u_1-x_2-u_2+y_2)^+}{\epsilon},1\big)\zeta=0 \quad \text{ on } \, (\Sigma\times Q)\cup
(Q\times \Sigma)$$
 and the functions $h(x_1)$ and $u_1$ (resp. $h(y_1)$ and $u_2$) do not depend on $y_2$  (resp.
 $x_2$), we have
\begin{eqnarray}\label{e3.1}
&&\int_{Q\times
Q}h(x_1)a(u_{1x_2})\big(\min\big(\frac{(u_1-x_2-u_2+y_2)^+}{\epsilon},1\big)\zeta\big)_{y_2}\,
dx dt dy ds=0,\\
 &&\int_{Q\times
Q}h(y_1)a(u_{2y_2})\big(\min\big(\frac{(u_1-x_2-u_2+y_2)^+}{\epsilon},1\big)\zeta\big)_{x_2}\,
dx dt dy ds=0,\\
&&u_{2x_2}=u_{1y_2}=0 \quad \text{ a.e. in }   Q\times Q.
\end{eqnarray}
Subtracting (3.5) from and (3.3) and using (3.6)-(3.8), we get
\begin{eqnarray}\label{e3.1}
&&\int_{Q\times Q}\Big\{\big[h(x_1)a((\partial_{x_2}+\partial_{y_2})u_{1})-h(y_1)a((\partial_{x_2}+\partial_{y_2})u_{2})\big]\nonumber\\
&&\times(\partial_{x_2}+\partial_{y_2})\big(\min\big(\frac{(u_1-x_2+y_2-u_2)^+}{\epsilon},1\big)\zeta\big)\nonumber\\
&&+(g_2-g_1)^+\big[(h(x_1)\zeta_{x_2}+h(y_1)\zeta_{y_2})a(1)-\zeta_{t}-\zeta_s\big]\Big\}\,
dx dtdyds\nonumber\\
&&\leq \int_{Q}C(u_1,u_2-y_2,\zeta)\, dy ds \nonumber
\end{eqnarray}
and taking into account the condition (1.3), it can be written as
\begin{eqnarray}\label{e3.1}
&&\int_{Q\times
Q}\Big\{\big[h(x_1)a((\partial_{x_2}+\partial_{y_2})u_{1})-h(y_1)
a((\partial_{x_2}+\partial_{y_2})u_{2})\big](\partial_{x_2}+\partial_{y_2})\zeta\nonumber\\
&&\quad\quad\times\min\big(\frac{(u_1-x_2+y_2-u_2)^+}{\epsilon},1\big)\nonumber\\
&&\quad\quad
+(g_2-g_1)^+\big[h(y_1)a(1)(\zeta_{x_2}+\zeta_{y_2})-\zeta_{t}-\zeta_s\big]\Big\}\,
dx dtdyds\nonumber\\
&&+\int_{Q\times
Q}(g_2-g_1)^+(h(x_1)-h(y_1))a(1)\zeta_{x_2}\, dx dt dy ds \nonumber\\
&&+\int_{Q\times
Q}(h(x_1)-h(y_1))a((\partial_{x_2}+\partial_{y_2})u_{2})\zeta  \nonumber\\
&&\quad\quad\times
(\partial_{x_2}+\partial_{y_2})\big(\min\big(\frac{(u_1-x_2+y_2-u_2)^+}{\epsilon},1\big)\big)\,dx dt dy ds\nonumber\\
&&+\int_{Q\times (Q\cap\{u_2=y_2\})}h(x_1)(a(u_{1x_2})-
a(1))\zeta\big(\min\big(\frac{u_1-x_2}{\epsilon},1\big)\big)_{x_2}\,dx dt dy ds\nonumber\\
&&\leq \int_{Q}C(u_1,u_2-y_2,\zeta)\, dy ds.
\end{eqnarray}
Let us consider the following change of variables:
\begin{eqnarray}\label{e3.1}
&&z=\frac{x+y}{2}, \quad\tau=\frac{t+s}{2},
\quad\sigma=\frac{x-y}{2}, \quad\theta=\frac{t-s}{2},
\end{eqnarray}
 let $J_1$ and
$J_2$ denote, respectively,  the domains of the variables
$z_2=\frac{x_2+y_2}{2}$ and  $\sigma_2=\frac{x_2-y_2}{2}$  and  let
$I$ denote  the domain of the variables $x_2$ and  $y_2$. Set
\begin{eqnarray}\label{e3.1}
&&Q_1=(A,B)\times J_1\times(0,T), \,
Q_2=\Big(\frac{A-B}{2},\frac{B-A}{2}\Big)\times
J_2\times\Big(-\frac{T}{2},\frac{T}{2}\Big), \nonumber\\
&&Q_3=(A,B)\times I\times(0,T) \, \text{ and }  \,
Q_4=\Big(\frac{A-B}{2},\frac{B-A}{2}\Big)\times
I\times\Big(-\frac{T}{2},\frac{T}{2}\Big).\nonumber
\end{eqnarray}
 Insert (3.10)  in (3.9) yields
\begin{eqnarray}\label{e3.1}
&&E_{\epsilon, \delta_2, \delta_1}+F_{\delta_2, \delta_1}+G_{\epsilon, \delta_2, \delta_1} \nonumber\\
&&+\int_{Q\times (Q\cap\{u_2=y_2\})}h(x_1)(a(u_{1x_2})-
a(1))\zeta\big(\min\big(\frac{u_1-x_2}{\epsilon},1\big)\big)_{x_2}\,dx dt dy ds\nonumber\\
&&\leq \int_{Q}C(u_1,u_2-y_2,\zeta)\, dy ds,
\end{eqnarray}
 where
 \begin{eqnarray}
 E_{\epsilon, \delta_2, \delta_1}&=&\int_{Q_1\times
Q_2}\Big\{\big[h(z_1+\sigma_1)a(\hat{u}_{1z_2})-h(z_1-\sigma_1)
a(\hat{u}_{2z_2})\big]\hat{\zeta}_{z_2}\min\big(\frac{(\hat{u}_1-\hat{u}_2-2\sigma_2)^+}{\epsilon},1\big)
\nonumber\\
&&\quad\quad+(\hat{g}_2-\hat{g}_1)^+\big(h(z_1-\sigma_1)a(1)\hat{\zeta}_{z_2}-\hat{\zeta}_{\tau}\big)\Big\}\,
dz d\tau d\sigma d\theta,\nonumber\\
F_{\delta_2,\delta_1}&=&\int_{Q_3\times
Q_4}(\bar{g}_2-\bar{g}_1)^+(h(z_1+\sigma_1)
-h(z_1-\sigma_1))a(1)\bar{\zeta}_{x_2}\,dz_1 dx_2d\tau d\sigma_1 dy_2 d\theta, \nonumber\\
G_{\epsilon, \delta_2, \delta_1}&=&\int_{Q_1\times
Q_2}(h(z_1+\sigma_1)-h(z_1-\sigma_1))a(\hat{u}_{2z_2})\hat{\zeta}
\big(\min\big(\frac{(\hat{u}_1-\hat{u}_2-2\sigma_2)^+}{\epsilon},1\big)\big)_{z_2}\,dz
d\tau d\sigma d\theta\nonumber
 \end{eqnarray}
 with
\begin{eqnarray}
&&\hat{u}_1=u_1(z+\sigma,\tau+\theta), \quad
\hat{u}_2=u_2(z-\sigma,\tau-\theta),
\quad\hat{\zeta}=\xi(z,\tau)\rho_{1,\delta_1}(\theta)\rho_{2,\delta_1}(\sigma_1)\rho_{3,\delta_2}(\sigma_2),\nonumber\\
&&\hat{g}_1=g_1(z+\sigma,\tau+\theta), \quad
\hat{g}_2=g_2(z-\sigma,\tau-\theta), \quad
\overline{u}_{1}=u_1(z_1+\sigma_1,x_2,\tau+\theta),\nonumber\\
&&\bar{g}_1=g_1(z_1+\sigma_1,x_2,\tau+\theta), \quad
\bar{g}_2=g_2(z_1-\sigma_1,y_2,\tau-\theta),\nonumber\\
&&\bar{\zeta}=\xi(z_1,\frac{x_2+y_2}{2},\tau)\rho_{1,\delta_1}(\theta)
\rho_{2,\delta_1}(\sigma_1)\rho_{3,\delta_2}(\frac{x_2-y_2}{2}).\nonumber
\end{eqnarray}
Since $h$ is a Lipschitz continuous function and
$\text{supp}(\rho_{2,\delta_1})\subset (-\delta_1,\delta_1)$, there
exists a constant $C$ such that
\begin{eqnarray}
|F_{\delta_2,\delta_1}|&\leq& 2C\int_{Q_3\times
Q_4}|\sigma_1|(\bar{g}_2-\bar{g}_1)^+a(1)|\bar{\zeta}_{x_2}|\,dz_1
dx_2d\tau d\sigma_1 dy_2 d\theta\nonumber\\
&\leq&2C\delta_1\int_{Q_3\times
Q_4}(\bar{g}_2-\bar{g}_1)^+a(1)|\bar{\zeta}_{x_2}|\,dz_1
dx_2d\tau d\sigma_1 dy_2 d\theta\nonumber\\
&:=&2C\delta_1W_{\delta_2,\delta_1}^1,\\
|G_{\epsilon, \delta_2, \delta_1}|&\leq&2C\delta_1\int_{Q_1\times
Q_2}
\Big|a(\hat{u}_{2z_2})\big(\min\big(\frac{(\hat{u}_1-\hat{u}_2-2\sigma_2)^+}{\epsilon},1\big)\big)_{z_2}\Big|\hat{\zeta}\,dz
d\tau d\sigma d\theta\nonumber\\
&:=&2C\delta_1W_{\epsilon, \delta_2, \delta_1}^2.
\end{eqnarray}
Notice that, $(W_{\delta_2,\delta_1}^1)_{\delta_1>0}$ and
$(W_{\epsilon, \delta_2, \delta_1}^2)_{\delta_1>0}$ are bounded,
then, passing to the limit in (3.12)-(3.13) as
$\delta_1\rightarrow0,$ we obtain
\begin{eqnarray}
\lim_{\delta_1\rightarrow0}(F_{\delta_2,\delta_1})=\lim_{\delta_1\rightarrow0}(G_{\epsilon,
\delta_2, \delta_1})=0.
\end{eqnarray}
On the other hand,
\begin{eqnarray}
 \lim_{\delta_2\rightarrow0}(\lim_{\delta_1\rightarrow0}(E_{\epsilon,\delta_2,\delta_1}))&=&\int_{Q\times
Q}\Big\{h(z_1)(a(u_{1z_2})-a(u_{2z_2}))\xi_{z_2}\min\big(\frac{(u_1-u_2)^+}{\epsilon},1\big)
\nonumber\\
&&\quad\quad+(g_2-g_1)^+(h(z_1)a(1)\xi_{z_2}-\xi_{\tau})\Big\}\, dz
d\tau,
 \end{eqnarray}
 where $u_1=u_1(z,\tau), \,  u_2=u_2(z,\tau), \, g_1=g_1(z,\tau), \,
g_2=g_2(z,\tau)$ and $\xi=\xi(z,\tau).$ Now, since, by  (2.9) and
the Lebesgue theorem, we have
\begin{eqnarray}\label{e3.1}
&&\int_{Q}C(u_1,u_2-y_2,\zeta)\, dy ds\nonumber\\
&&=\lim_{\epsilon\rightarrow0}\int_{Q\cap\{u_2=y_2\}}\Big\{\int_{Q}h(x_1)(a(u_{1x_2})-
a(1))\zeta\big(\min\big(\frac{u_1-x_2}{\epsilon},1\big)\big)_{x_2}\,dx
dt\Big\}dy ds,\nonumber
\end{eqnarray}
we obtain by letting successively $\delta_1\rightarrow0$,
$\delta_2\rightarrow0$, $\epsilon\rightarrow0$ in (3.11) and using
(3.14)-(3.15),
\begin{eqnarray}\label{e3.1}
&&\int_{Q}\big\{\chi_{\{u_1-u_2\geq0\}}h(z_1)(a(u_{1z_2})-
a(u_{2z_2}))\xi_{z_2}+(g_2-g_1)^+(h(z_1)a(1)\xi_{z_2}-\xi_{\tau})\big\}\,
dz d\tau\leq 0,\nonumber
\end{eqnarray}
where $\chi_{\{u_1-u_2\geq0\}}$ denotes the characteristic function
of the set $\{u_1-u_2\geq0\}.$ This leads to (3.1) for $i=1$. If one
exchanges the roles of $(u_1,g_1)$ and $(u_2,g_2)$, one  also
obtains (3.1) for $i=2$.  \qed

\begin{lemma}
Let $B$ be a bounded open subset of $\mathbb{R}^{2}$ such that
either $B\cap\Gamma=\emptyset$  or $B\cap\Gamma$  is a Lipschitz
graph. Let $(u_1,g_1)$ and $(u_2,g_2)$ be two solutions of $(P)$ and
let $\overline{g}$ be a function of $L^{\infty}(Q)$ such that
\begin{eqnarray}\label{e3.1}
&&0\leq \overline{g}\leq g_1, g_2 \text{ a.e. in }  Q, \quad
h(x_1)a(1)\overline{g}_{x_2}-\overline{g}_{t}=0 \text{ in }
\mathcal{D}^{\prime}(Q).
\end{eqnarray}
Then, for all $\xi\in \mathcal{D}(B\times(0,T)), \, \xi\geq0, \,
\text{supp}(\xi)\cap(\sigma_1\cup \Sigma_4)=0$ and for $i,j=1,2, \,
i\neq j,$ we have
\begin{eqnarray}\label{e3.1}
&&\int_{Q}\big\{h(x_1)(a(u_{ix_2})-a(u_{mx_2})-(g_j-\overline{g})^+a(1))\xi_{x_2}\nonumber\\
&&\quad\quad+(g_j-\overline{g})^+\xi_t\big\}\,dxdt\leq0.
\end{eqnarray}
\end{lemma}
\n \emph{Proof.} Let $(u_1,g_1)$ and $(u_2,g_2)$ be two solutions of
 $(P)$ and let $\xi$ be the function defined in Lemma 3.2. We
 define
 \begin{eqnarray}\label{e3.1}
&&\forall (x,t,y,s)\in \overline{Q\times Q}:\nonumber\\
&&\zeta(x,t,y,s)=\xi(\frac{x+y}{2},\frac{t+s}{2})
\rho_{1,\delta_1}(\frac{t-s}{2})\rho_{2,\delta_2}(\frac{x_1-y_1}{2})\rho_{3,\delta_3}(\frac{x_2-y_2}{2}),\nonumber
\end{eqnarray}
where $\delta_1$, $\delta_2$, $\delta_3$  are  positive real
numbers, $\rho_{1,\delta_1}, \,\rho_{2,\delta_2},
\,\rho_{3,\delta_3}\in \mathcal{D}(\mathbb{R}),\, \rho_{1,\delta_1},
\, \,\rho_{2,\delta_2}, \, \rho_{3,\delta_3}\geq0$ in $\mathbb{R},\,
\int_{\mathbb{R}}\rho_{1,\delta_1}(t)\,dt=\int_{\mathbb{R}}\rho_{2,\delta_2}(t)\,dt=\int_{\mathbb{R}}\rho_{3,\delta_3}(t)\,dt=1,
\, \text{supp}(\rho_{1,\delta_1})\subset(-\delta_1,\delta_1),\,
\text{supp}(\rho_{2,\delta_2})\subset(-\delta_2,\delta_2),
\,\text{supp}(\rho_{3,\delta_3})\subset (-\delta_3,\delta_3)$  and
\begin{eqnarray}\label{e3.1}
&&\forall (x,y)\in (B\cap\Omega)\times (B\backslash\Omega), \quad
\rho_{2,\delta_2}(\frac{x_1-y_1}{2})\rho_{3,\delta_3}(\frac{x_2-y_2}{2})=0.\nonumber
\end{eqnarray}
For  $\delta_1$, $\delta_2$ and $\delta_3$ small enough, we have
$\zeta\in \mathcal{D}(B\times(0,T)\times B\times(0,T))$ and
\begin{eqnarray}\label{e3.1}
&&\zeta=0 \quad \text{ on } \, (\Sigma\times Q)\cup (Q\times
(\sigma_1\cup\Sigma_4)).\nonumber
\end{eqnarray}
On the other hand, since $\text{supp}(\xi)\cap
(\sigma_1\cap\Sigma_4)=\emptyset$ and if we suppose that
$\text{supp}(\xi)\cap \Sigma_3\neq\emptyset,$ we can find $r_0\in
\Big(0,\displaystyle{\min_{\text{supp}(\xi)\cap
\Sigma_3}\varphi}\Big)$ and $\Psi\in C^{\infty}(\mathbb{R})\cap
C^{0,1}(\mathbb{R})$ such that $\Psi^{\prime}\geq0,\, \Psi(r)=0$ if
$r\leq0$ and  $\Psi(r)=1$ if $r\geq r_0,$ and  this function $\Psi$
satisfies for $\delta_1$, $\delta_2$ and $\delta_3$ small enough,
\begin{eqnarray}\label{e3.1}
&&(1-\Psi(u_2-y_2))\zeta=0 \quad \text{ on } \, (\Sigma\times Q)\cup
(Q\times \Sigma_2).
\end{eqnarray}
If $\text{supp}(\xi)\cap \Sigma_3=\emptyset,$ we choose $\Psi=0.$
Now, applying Lemma 2.2 to $(u_1,g_1)$ with $k=u_2(y,s)-y_2, \,
\lambda=g_2(y,s), \,\xi(x,t)=\xi_2(x,t)=\zeta(x,t,y,s),
\,\xi_1(x,t)=0$ and $\chi=\overline{g},$ we obtain for a.e.
$(y,s)\in Q,$
\begin{eqnarray}\label{e3.1}
&&\int_{Q}\Big\{h(x_1)a(u_{1x_2})\big(\min\big(\frac{(u_1-x_2+y_2-u_2)^+}{\epsilon},1\big)\zeta\big)_{x_2}\nonumber\\
&&+(g_2-\overline{g})^+(h(x_1)a(1)\zeta_{x_2}-\zeta_{t})\Big\}\, dx
dt\leq 0\nonumber
\end{eqnarray}
and integrating over $Q$, we get
\begin{eqnarray}\label{e3.1}
&&\int_{Q\times Q}\Big\{h(x_1)a(u_{1x_2})\big(\min\big(\frac{(u_1-x_2+y_2-u_2)^+}{\epsilon},1\big)\zeta\big)_{x_2}\nonumber\\
&&+(g_2-\overline{g})^+(h(x_1)a(1)\zeta_{x_2}-\zeta_{t})\Big\}\, dx
dtdyds\leq 0.
\end{eqnarray}
Similarly, for a.e. $(x,t)\in Q,$ we apply Lemma 2.3 to $(u_2,g_2)$
with $k=u_1(x,t)-x_2, \, \lambda=\overline{g}(x,t),
\,\xi(y,s)=\zeta(x,t,y,s),$ then, we integrate over $Q$ to obtain
\begin{eqnarray}\label{e3.1}
&&-\int_{Q\times
Q}\Big\{h(y_1)\big[a(u_{2y_2})-\overline{g}a(1)\big]
\big(\min\big(\frac{(u_1-x_2+y_2-u_2)^+}{\epsilon},1\big)(1-\Psi(u_2-y_2))\zeta\big)_{y_2}\nonumber\\
&&\qquad\qquad
-(g_2-\overline{g})^+(h(y_1)a(1)\zeta_{y_2}-\zeta_{s})\Big\}\, dx
dtdyds\leq 0.
\end{eqnarray}
On the other hand, by Corollary 2.1,  we have
\begin{eqnarray}\label{e3.1}
&&\int_{Q\times
Q}h(x_1)a(u_{1x_2})\big(\min\big(\frac{(u_1-x_2+y_2-u_2)^+}{\epsilon},1\big)\zeta\big)_{x_2}\,
dx dtdyds= 0
\end{eqnarray}
and the use of (3.16) and (3.18) leads to
\begin{eqnarray}\label{e3.1}
\int_{Q\times
Q}\overline{g}h(y_1)a(1)\big(\min\big(\frac{(u_1-x_2+y_2-u_2)^+}{\epsilon},1\big)(1-\Psi(u_2-y_2))\zeta\big)_{x_2}\,
dx dtdyds= 0.\nonumber\\
\end{eqnarray}
Addition of (3.19), (3.20), (3.21) and (3.22) yields
\begin{eqnarray}\label{e3.1}
&&K_{\epsilon,\delta_1,\delta_3,\delta_2}+L_{\delta_1,
\delta_3,\delta_2}+M_{\epsilon,\delta_1,\delta_3,\delta_2}\leq0,
\end{eqnarray}
where
\begin{eqnarray}\label{e3.1}
K_{\epsilon,\delta_1,\delta_3,\delta_2}&=&\int_{Q\times
Q}\Big\{\big[h(x_1)a(u_{1x_2})-h(y_1)a(u_{2y_2})\big]
(\partial_{x_2}+\partial_{y_2})\big(\min\big(\frac{(u_1-x_2+y_2-u_2)^+}{\epsilon},1\big)\zeta\big)\nonumber\\
&&+(g_2-\overline{g})^+(\zeta_{x_2}+\zeta_{y_2})h(y_1)a(1)-\zeta_{t}-\zeta_s\big)\Big\}\,
dx dt dy d s, \nonumber\\
L_{\delta_1, \delta_3, \delta_2}&=&\int_{Q\times
Q}(g_2-\overline{g})^+(h(x_1)-h(y_1))a(1)\zeta_{x_2}\, dx dt
dy ds, \nonumber\\
M_{\epsilon,\delta_1,\delta_3,\delta_2}&=&\int_{Q\times
Q}h(y_1)\big[a((\partial_{x_2}+\partial_{y_2})u_{2})-\overline{g}a(1)\big]\nonumber\\
&&\qquad\qquad\times(\partial_{x_2}+\partial_{y_2})\big(\min\big(\frac{(u_1-x_2+y_2-u_2)^+}{\epsilon},1\big)\Psi(u_2-y_2)\zeta\big)\,dx
dt dy d s, \nonumber\\
&&-\int_{Q\times
Q}\big[h(x_1)a((\partial_{x_2}+\partial_{y_2})u_{1})-h(y_1)\overline{g}a(1)\big]\nonumber\\
&&\qquad\qquad\times(\partial_{x_2}+\partial_{y_2})\big(\min\big(\frac{(u_1-x_2+y_2-u_2)^+}{\epsilon},1\big)\zeta\big)\,dx
dt dy d s. \nonumber
\end{eqnarray}
Passing to the limit in $L_{\delta_1, \delta_3,  \delta_2}$ as
$\delta_2\rightarrow0,$ we arrive at
\begin{eqnarray}\label{e3.1}
&&\lim_{\delta_2\rightarrow 0}(L_{\delta_1, \delta_3, \delta_2})=0.
\end{eqnarray}
On the other hand,
\begin{eqnarray}\label{e3.1}
&&\lim_{\delta_3\rightarrow 0}(\lim_{\delta_2\rightarrow
0}(M_{\epsilon,\delta_1,\delta_3,\delta_2}))\nonumber\\
 &&=\int_0^T\int_{Q}h(x_1)\big[a(u_{2x_2})-\overline{g}a(1)\big]\nonumber\\
&&\qquad\times\big(\min\big(\frac{(u_1-u_2)^+}{\epsilon},1\big)\Psi(u_2-x_2)\xi\big)_{x_2}\rho_{1,\delta_1}\,dx
dt d s \nonumber\\
&&-\int_0^T\int_{Q}h(x_1)\big[a(u_{1x_2})-\overline{g}a(1)\big]\nonumber\\
&&\qquad\times\big(\min\big(\frac{(u_1-u_2)^+}{\epsilon},1\big)\xi\big)_{x_2}\rho_{1,\delta_1}\,dx
dt d s:=S_{\epsilon, \delta_1},
\end{eqnarray}
where $u_1=u_1(x,t)$, $u_2=u_2(x,s)$
$\overline{g}=\overline{g}(x,t)$, $\xi=\xi(x,\frac{t+s}{2})$ and
$\rho_{1,\delta_1}=\rho_{1,\delta_1}(\frac{t-s}{2}).$  Applying
Lemma 2.1  to
$F(z_1,z_2)=\min(\frac{(z_1-z_2)^+}{\epsilon},1)(1-\Psi(z_2))$ with
$v=u_2-x_2$ and taking into account
$$(1-\Psi(u_2(x,s)-x_2))\xi(x,\frac{t+s}{2})\rho_{1,\delta_1}(\frac{t-s}{2})=0$$
for all $(x,t,s)\in \Sigma_2\times(0,T),$ we get
\begin{eqnarray}\label{e3.1}
&&\int_0^T\int_{Q}h(x_1)\big[a(u_{1x_2})-g_1a(1)\big]\nonumber\\
&&\times\big(\min\big(\frac{(u_1-u_2)^+}{\epsilon},1\big)(1-\Psi(u_2-x_2))\xi\big)_{x_2}\rho_{1,\delta_1}\,dx
dt d s=0.
\end{eqnarray}
Using (3.26) and the fact that $\overline{g}\leq g_1$ and
$\Psi(0)=0,$ we obtain from (3.25),
\begin{eqnarray}\label{e3.1}
 S_{\epsilon, \delta_1}&=&\int_0^T\int_{Q}h(x_1)\big[a(u_{2x_2})-g_2a(1)\big]\nonumber\\
&&\times\big(\min\big(\frac{(u_1-u_2)^+}{\epsilon},1\big)\Psi(u_2-x_2)\xi\big)_{x_2}\rho_{1,\delta_1}\,dx
dt d s \nonumber\\
&&-\int_0^T\int_{Q}h(x_1)\big[a(u_{1x_2})-g_1a(1)\big]\nonumber\\
&&\times\big(\min\big(\frac{(u_1-u_2)^+}{\epsilon},1\big)\Psi(u_2-x_2)\xi\big)_{x_2}\rho_{1,\delta_1}\,dx
dt d s.
\end{eqnarray}
Notice that
$$\pm\big(\min(\frac{(u_1-\phi(x,s))^+}{\epsilon},1)-\min(\frac{(\phi(x,t)-\phi(x,s))^+}{\epsilon},1)\big)\Psi(\phi(x,s)-x_2)\xi\rho_{1,\delta_1}$$
are test functions for $(P)$ corresponding to $(u_2,g_2)$.  In
addition, applying Lemma 2.1 to $u_2$ with $v=u_1-x_2$,
$F(z_1,z_2)=\min(\frac{(z_2-z_1)^+}{\epsilon},1)(1-\Psi(z_1))$ and
$F(z_1,z_2)=\min(\frac{(z_2-z_1)^+}{\epsilon},1)$, subtracting one
equation from the other and  taking into account $\Psi(0)=0$,
$g_2(u_2-x_2)=0$ a.e. in $Q$, we deduce that
\begin{eqnarray}\label{e3.1}
&&\int_0^T\int_{Q}h(x_1)\big[a(u_{2x_2})-g_2a(1)\big]\nonumber\\
&&\times\big(\min\big(\frac{(u_1-u_2)^+}{\epsilon},1\big)\Psi(u_2-x_2)\xi\big)_{x_2}\rho_{1,\delta_1}\,dx
dt d s \nonumber\\
&&=\int_0^T\int_{Q}h(x_1)\Big\{\big[a(u_{2x_2})-g_2a(1)\big]\nonumber\\
&&\times\big(\min\big(\frac{(\phi(x,t)-\phi(x,s))^+}{\epsilon},1\big)\Psi(\phi(x,s)-x_2)\xi\big)_{x_2}\rho_{1,\delta_1}\nonumber\\
&&+g_2\big(\min\big(\frac{(\phi(x,t)-\phi(x,s))^+}{\epsilon},1\big)\Psi(\phi(x,s)-x_2)\xi\rho_{1,\delta_1}\big)_{s}\Big\}\,dx
dt d s.
\end{eqnarray}
Similarly, if we apply Lemma 2.1 to $u_1$ with $v=u_2-x_2$ and
$F(z_1,z_2)=\min(\frac{(z_1-z_2)^+}{\epsilon},1)\Psi(z_2)$, we get
\begin{eqnarray}\label{e3.1}
&&\int_0^T\int_{Q}h(x_1)\big[a(u_{1x_2})-g_1a(1)\big]\nonumber\\
&&\times\big(\min\big(\frac{(u_1-u_2)^+}{\epsilon},1\big)\Psi(u_2-x_2)\xi\big)_{x_2}\rho_{1,\delta_1}\,dx
dt d s \nonumber\\
&&=\int_0^T\int_{Q}h(x_1)\Big\{\big[a(u_{1x_2})-g_1a(1)\big]\nonumber\\
&&\times\big(\min\big(\frac{(\phi(x,t)-\phi(x,s))^+}{\epsilon},1\big)\Psi(\phi(x,s)-x_2)\xi\big)_{x_2}\rho_{1,\delta_1}\nonumber\\
&&+g_1\big(\min\big(\frac{(\phi(x,t)-\phi(x,s))^+}{\epsilon},1\big)\Psi(\phi(x,s)-x_2)\xi\rho_{1,\delta_1}\big)_{t}\Big\}\,dx
dt d s.
\end{eqnarray}
Since $g_1$ (resp. $g_2$) does not depend on $s$ (resp. $t$), we
obtain from (3.27), (3.28) and (3.29),
\begin{eqnarray}\label{e3.1}
 S_{\epsilon, \delta_1}&=&\int_0^T\int_0^T\int_{Q}\Big\{h(x_1)\big[a(u_{2x_2})-a(u_{1x_2})+(g_1-g_2)a(1)\big]\nonumber\\
&&\times\big(\min\big(\frac{(\phi(x,t)-\phi(x,s))^+}{\epsilon},1\big)\Psi(\phi(x,s)-x_2)\xi\big)_{x_2} \nonumber\\
&&-(g_1-g_2)(\partial_t+\partial_s)\big(\min\big(\frac{(\phi(x,t)-\phi(x,s))^+}{\epsilon},1\big)\Psi(\phi(x,s)-x_2)\xi\big)\Big\}\rho_{1,\delta_1}\,dx
dt d s.\nonumber\\
\end{eqnarray}
We may then pass to the limit in (3.30) as $\delta_1\rightarrow0$ to
deduce
\begin{eqnarray}
\lim_{\delta_1\rightarrow0}(\lim_{\delta_3\rightarrow 0}(\lim_{\delta_2\rightarrow
0}(M_{\epsilon,\delta_1,\delta_3,\delta_2}))=\lim_{\delta_1\rightarrow0}(S_{\epsilon,
\delta_1})=0.
\end{eqnarray}
 Thus, for $K_{\epsilon,\delta_1,\delta_3,\delta_2}$, we use (1.3) and pass successively to the
 limit as $\delta_2\rightarrow0$, $\delta_3\rightarrow0$, $\delta_1\rightarrow0$,
 $\epsilon\rightarrow0$ to get
\begin{eqnarray}\label{e3.1}
&&\int_{Q\times
Q}\Big\{\chi_{\{u_1-u_2\geq0\}}h(x_1)\big[a(u_{1x_2})-a(u_{2x_2})\big]\zeta_{x_2}\nonumber\\
&&+(g_2-\overline{g})^+\big(h(x_1)a(1)\zeta_{x_2}-\zeta_{t})\Big\}\,
dx dt\leq
\liminf_{\epsilon\rightarrow0}(\lim_{\delta_1\rightarrow0}(\lim_{\delta_3\rightarrow
0}(\lim_{\delta_2\rightarrow
0}(K_{\epsilon,\delta_1,\delta_3,\delta_2})))).\nonumber\\
\end{eqnarray}
Now, by letting successively $\delta_2\rightarrow0$,
$\delta_3\rightarrow0$, $\delta_1\rightarrow0$,
$\epsilon\rightarrow0$ in (3.23) and
 using (3.24), (3.31) and (3.32), we obtain (3.17) for $i=1$ and $j=2$. Since we can exchange the roles of $(u_1,g_1)$ and $(u_2,g_2)$,
 we can also obtain (3.17) for $i=2$ and $j=1$. \qed

\n \emph{Proof of Theorem 3.1.} Using Lemmas 3.1, 3.2 and applying
arguments similar to [\cite{[Ly1]}, Lemmas 5.6, 5.7 and Theorem
5.8], we arrive at the following comparison of solutions,
\begin{eqnarray}\label{e3.1}
&&\int_{Q}\big\{h(x_1)(a(u_{ix_2})-a(u_{mx_2})-(g_i-g_M)a(1))\xi_{x_2}\nonumber\\
&&\quad\quad+(g_i-g_M)\xi_t\big\}\,dxdt\leq0,\\
&&i=1,2, \, \forall \xi \in \mathcal{D}(B\times(0,T)), \, \xi\geq0,
\, \xi(x,0)=\xi(x,T)=0 \text{ a.e. in } \Omega.\nonumber
\end{eqnarray}
If we choose $\xi \in \mathcal{D}(0,T), \, \xi\geq0$ in (3.33), we
get
$$\int_{Q}(g_i-g_M)\xi_t\, dxdt\leq0,$$
which can be written as
$$\frac{d}{dt}\int_{\Omega}(g_M-g_i)\, dx\leq0 \quad \text{ in } \mathcal{D}^{\prime}(0,T).$$
Since $g_i\in C^0([0,T];L^1(\Omega))$ (see \cite{[Ly2]}) and
$g_1(x,0)=g_2(x,0)=g_0(x)$ a.e. $x\in \Omega$, we obtain
$$\int_{\Omega}(g_M-g_i)\, dx=0 \quad \text{ in } [0,T],$$
which leads to
\begin{eqnarray}\label{e3.1}
g_1=g_2=g_M \quad \text{ a.e. in } Q.
\end{eqnarray}
Insert (3.34) in (3.33) yields
\begin{eqnarray}\label{e3.1}
\forall \xi \in \mathcal{D}(\overline{\Omega}\times(0,T)), \,
\xi\geq0: \quad
\int_{Q}h(x_1)(a(u_{ix_2})-a(u_{mx_2}))\xi_{x_2}\,dxdt\leq0.
\end{eqnarray}
From (1.3)-(1.4) and the fact that (3.35) remains true for
$\xi=u_i-u_m$, we deduce that $(u_i-u_m)_{x_2}=0$ a.e. in $Q$. Since
$u_i-u_m=0$ on $\Sigma_2$, we can extend $u_i-u_m$ to
$\mathbb{R}\times(D,+\infty)\times(0,T)$ by $0$ and still denote by
$u_i-u_m.$ Thus, for a.e. $(x_1,t)\in \mathbb{R}\times(0,T)$, there
exists  $w\in C^0([D,+\infty))$ such that
$w(x_2)=(u_i-u_m)(x_1,x_2,t)$ a.e. $x_2\in (D,+\infty)$ and
$$\forall z_1, z_2\in [D,+\infty): \quad  w(z_1)-w(z_2)=\int_{z_2}^{z_1}(u_i-u_m)_{x_2}(x_1,z,t)\,dz=0,$$
which means that $w=c$ in $[D,+\infty)$ for some constant $c\geq0.$
Due to $w(x_2)=0$ for $x_2$ large enough,  it follows that $w=0$ in
$[D,+\infty)$, and hence $u_i=u_m$ a.e. in $Q$ for $i=1,2.$  Thus
the proof is complete. \qed

\begin{remark}\label{r2.1}
 If $E$, $F$ are real numbers such that $F>E$,  $n=3$ and
 $\Gamma_1=[A,B]\times[E,F]$, the obtained uniqueness result remains true if we
 replace $x_1$ by $x^\prime=(x_1,x_2)\in [A,B]\times[E,F]$ and $x_2$ by $x_3$, where
 $x=(x^\prime,x_3)=(x_1,x_2,x_3)$ is a generic point of $\Omega\subset
 \mathbb{R}^3.$ In this way, our uniqueness result can be also extended
 to a dam of $\mathbb{R}^n$ with $n\geq4.$
\end{remark}

\end{document}